\documentstyle[12pt,righttag,ctagsplt]{amsart}
    \newtheorem{prop}{Proposition}
    \newtheorem{th}[prop]{Theorem}
    \newtheorem{cor}[prop]{Corollary}
    \newtheorem{lem}[prop]{Lemma}

    \theoremstyle{definition}
     
    \theoremstyle{definition}
    
    \theoremstyle{remark}
    \newtheorem{rem}{Remark} 
    \newcommand{\bbN}{{\Bbb{N}}}
    \newcommand{\bbR}{{\Bbb{R}}}
    
    \newcommand{\calA}{{\cal{A}}}

    \newcommand{\calL}{{\cal{L}}}
    
    \newcommand{\al}{\alpha}
    \newcommand{\be}{\beta}

    \newcommand{\e}{{\varepsilon}}

    \renewcommand{\span}{\operatorname{span}}
    
    \newcommand{\rank}{\operatorname{rank}}
    \newcommand{\card}{\operatorname{card}}

    \newcommand{\sgn}{\operatorname{sgn}}

    \newcommand{\disp}{\displaystyle}
    \newcommand{\lb}{\label}
    \newcommand{\ti}{\tilde}
    \newcommand{\wt}{\widetilde}
    \newcommand{\emp}{\emptyset}

    \newcommand{\lra}{\longrightarrow}
    
\newcommand{\starr}{={\hspace{-5mm}}\raisebox{-3mm}{\small (*)}}
\newcommand{\starrr}{={\hspace{-7mm}}\raisebox{-3mm}{\small (**)}}

    \newcommand{\DEF}{\buildrel {\mbox{\rm\small def}}\over =}

    \makeatletter
    \def\@currentlabel{2.1}\label{e:dispaa}
    \def\@currentlabel{2.21}\label{e:dispau}
    \def\@currentlabel{2.22}\label{e:dispav}
    \def\@currentlabel{2.23}\label{e:dispaw}
    \def\@currentlabel{2.24}\label{e:dispax}
    \makeatother
    \makeatletter
    \def\alphenumi{%
    \def\theenumi{\alph{enumi}}%
    \def\p@enumi{\theenumi}%
    \def\labelenumi{(\@alph\c@enumi)}}
    \makeatother
    \font\tenex=cmex10
    \newskip\ttglue
    \def\eightpoint{\def\rm{\fam0\eightrm}
    \textfont0=\eightrm \scriptfont0=\sixrm
    \scriptscriptfont0=\fiverm
    \textfont1=\eighti  \scriptfont1=\sixi
    \scriptscriptfont1=\fivei
    \textfont2=\eightsy  \scriptfont2=\sixsy
    \scriptscriptfont2=\fivesy
    \textfont3=\tenex  \scriptfont3=\tenex
    \scriptscriptfont3=\tenex
    \textfont\itfam=\eightit  \def\it{\fam\itfam\eightit}
    \textfont\slfam=\eightsl  \def\sl{\fam\slfam\eightsl}
    \textfont\ttfam=\eighttt  \def\tt{\fam\ttfam\eighttt}
    \textfont\bffam=\eightbf  \scriptfont\bffam=\sixbf
    \scriptscriptfont\bffam=\fivebf
    \def\bf{\fam\bffam\eightbf}
    \tt  \ttglue=.5em plus.25em minus.15em
    \normalbaselineskip=9pt
    \setbox\strutbox=\hbox{\vrule height7pt depth2pt width0pt}
    \let\sc=\sixrm  \let\big=\eightbig \normalbaselines\rm}
    \font\eightrm=cmr8 \font\sixrm=cmr6 \font\fiverm=cmr5
    \font\eighti=cmmi8  \font\sixi=cmmi6   \font\fivei=cmmi5
    \font\eightsy=cmsy8  \font\sixsy=cmsy6 \font\fivesy=cmsy5
    \font\eightit=cmti8  \font\eightsl=cmsl8
    \font\eighttt=cmtt8
    \font\eightbf=cmbx8  \font\sixbf=cmbx6 \font\fivebf=cmbx5
    \def\eightbig#1{{\hbox{$\textfont0=\ninerm\textfont2=\ninesy
    \left#1\vbox to6.5pt{}\right.\enspace$}}}

\begin{document}
\title{On isometric stability of complemented subspaces of ${\bold L_p}$}
\author{Beata Randrianantoanina}
\address{Department of Mathematics and Statistics
         \\ Miami University \\ Oxford, OH 45056}
\email{randrib@@muohio.edu}


\begin{abstract}
We show that Rudin-Plotkin isometry extension theorem in $L_p$ implies
that when $X$ and $Y$ are isometric subspaces of $L_p$ and $p$ is not an
even integer, $1 \leq p < \infty$, then $X$ is complemented in $L_p$ if and
only if $Y$ is; moreover the constants of complementation of $X$ and $Y$
are equal.  We provide examples demonstrating that this fact fails when $p$
is an even integer larger than 2.
\end{abstract} \maketitle

\section{Introduction}

Complemented subspaces of Banach spaces and in particular of $L_p$
have been intensely studied since the introduction of the notion of a Banach
space. Rosenthal \cite{R70} and others (see the survey \cite{force}) 
demonstrated that a class of
complemented subspaces of   $L_p$ is isomorphically very rich.
Moreover it was shown that complemented subspaces of   
$L_p$ are not  ``isomorphically stable'', i.e. there exist pairs of
isomorphic subspaces $X$, $Y$ of $L_p$ so that $X$ is complemented and
$Y$ is not. It is even possible to take $X$, $Y$ isomorphic to $L_p$
(see \cite{LR69a,R70,BDGJN77,B81}).

In this paper we observe that complemented subspaces of $L_p$ are
``isometrically stable'' when $p$ is not an even integer.  That is, we show
that if $X, \ Y$ are isometric subspaces of $L_p$ ($p \neq 4, \ 6, \
\dots$) and $X$ is complemented in $L_p$ then  $Y$ is also complemented in
$L_p$, moreover, the constant of complementation does not change.

This fact is a consequence of the Rudin-Plotkin isometry extension theorem
in $L_p$, which fails when $p$ is an even integer.  We show that also the
``isometric stability'' of complemented subspaces fails in $L_p$ when $p =
4, \ 6, \ \dots$, i.e. we construct isometric subspaces $U_p, \ V_p$ in
$L_p$ so that $V_p$ is complemented and $U_p$ is not complemented in $L_p$.

Our construction is based on the Rosenthal's construction of uncomplemented
copy of $X_p$ in $L_p$ \cite{R70} and thus we answer a question of
Rosenthal formulated in \cite[Remark~1 after Proposition~5]{R70} whether or
not $\widetilde X_{p,w}$ is isometric to a complemented subspace of $L_p$
(no if $p$ is not an even integer, yes if $p = 4$).

Our study is motivated by   results of Dor \cite{D75}, Schechtman
\cite{Sch79} and Alspach \cite{A83}, who showed that if $X$ is
$1$-complemented in $L_p, \ 1 \leq p < \infty$ and $Y$ is ($1 +
\e$)-isomorphic with $X$ (for $\e$ small enough,
depending only on $p$) then $Y$ is complemented in $L_p$ and the constant
of complementation of $Y$ tends to $1$ as $\e$ tends to $0$.  We
do not know whether such an almost-isometric stability is valid for any
class of complemented but not $1$-complemented subspaces of $L_p$.

We will use standard notation and facts from Banach space theory as may be
found in \cite{LT2}. The construction of the main example in 
Section~\ref{example}
relies heavily on Rosenthal's construction of the space $X_p$ \cite{R70}.

\medskip

{\bf{Acknowledgements:}}  I would like to thank Professors G. Schechtman
and R. Pol for many valuable discussions.

\section{Isometric stability of complemented subspaces of $L_p$ when $p$ is
not an even integer}

We start with recalling Rudin-Plotkin isometry extension theorem \cite{Rd76,
P70,P74} (cf. also \cite{H81,L78}).

\begin{th} \lb{RP}
Let $p \geq 1$, $p$ not an even integer, and let $H \subset L_p [0, 1]$.
If $T:H \lra L_p [0, 1]$ is a linear isometry, then there exists a linear
isometry $\wt T:L_p ([0, 1], \ \calA_0) \lra L_p [0, 1]$ such that $T =
\wt T|_{H}$.  Here $\calA_0$ is the smallest $\sigma$-subalgebra making
all functions from the space $H$ measurable. \end{th}

As an almost immediate consequence, we obtain:
\begin{cor} \lb{isometry}
Let $p \geq 1, \ p$ not an even integer, and let $X$ be a $K$-complemented
subspace of $L_p [0, 1]$.  Suppose that $Y \subset L_p [0, 1]$ is isometric
with $X$.  Then $Y$ is $K$-complemented in $L_p [0, 1]$. \end{cor}

\begin{pf}  Let $T:Y \lra X$ denote the isometry between $X$ and $Y$.  Then
by Theorem~\ref{RP}, there exists an extension of $T$, $$\wt T:L_p
([0,1], \ \calA_0) \lra L_p[0,1]$$ so that $\| \wt T \| = 1$ and $\wt T |_Y =
T$, where $\calA_0$ is the smallest $\sigma$-subalgebra, making all functions
from $Y$ measurable.  By \cite{Ando}, there exists a contractive projection
$P_1:L_p[0,1] \lra L_p ([0,1], \ \calA_0)$.  Let $P$ denote the projection
$P:L_p[0,1] \lra X$.

 Define $Q:L_p[0,1] \lra Y$ by
\begin{align*} Q &= T^{-1} P \wt T P_1.\\
\intertext{Then, for every $y \in Y$,}
Qy &= T^{-1} P \wt T P_1 y = T^{-1} P \wt T y = T^{-1} P (Ty) = T^{-1} (Ty)
= y
\end{align*}
so $Q$ is a projection, and $\| Q \| \leq \| T^{-1} \| \cdot \| P \| \cdot
\| \wt T \| \cdot \| P_1 \| = \| P \|$.
\end{pf}

Theorem~\ref{RP} is not true when $p$ is an even integer $(p \neq 2)$, and
also, Corollary~\ref{isometry} fails for $p = 4, \ 6, \ \dots$.  The next
section is devoted to the construction of the example.

\section{Example} \lb{example}

In this section, we will show that if $p$ is an even integer greater than
$2$ then Theorem~\ref{isometry} fails in $L_p$, i.e. we will show that
there exist subspaces $V_p, \ U_p \subset L_p$ such that $V_p$ and $U_p$
are isometric, $V_p$ is complemented in $L_p$, and $U_p$ is not
complemented in $L_p$.

Our construction is somewhat long and will be divided into several lemmas.

First we observe that when $p$ is an even integer, say $p = 2k$, then the
$p$-norm of the sum of independent symmetric random variables is determined
by the norms of the summands.  Namely we have:

\begin{lem} \lb{sum}
Suppose that $p = 2k, \ \{f_j\}^n_{j=1}$ are independent symmetric random
variables and let $f = \sum^n_{j=1} f_j$.  Then $\|f\|_p$ depends only on
the values of ${(\|f_j\|_{2m})^k_{m=1,}}^n_{j=1}$.

  Namely we have

$${\|f\|_p}^p = \sum_{k_1 + \ldots + k_n = k} {2k \choose 2k_1}  {2(k-k_1)
\choose 2k_2}  \ldots  {2(k-k_1-\ldots-k_{n-1}) \choose 2k_n} \prod^n_{j=1}
{\|f_j\|_{2k_j}}^{2k_j}.$$
\end{lem}

\begin{pf} 
\begin{align*}\lb{sum1}
{\|f\|_p}^p &= \int^1_0 \bigg(\sum^n_{j=1} f_j\bigg)^{2k} \\
&= \int^1_0 \sum_{k_1+\ldots+k_n=2k} {2k \choose k_1} \ {2k-k_1 \choose k_2}
\ldots  {k_n \choose k_n}  \cdot  {f_1}^{k_1}  \cdot {f_2}^{k_2}
\cdots {f_n}^{k_n} \\
&{\starr} \sum_{k_1+\ldots+k_n=2k} {2k \choose k_1} \ldots {k_n
\choose k_n} \prod^n_{j=1} \big(\int^1_0 {f_j}^{k_j}\big) \\
&{\starrr} \sum_{k_1+\ldots+k_n=k} {2k \choose 2k_1} \ldots {2k_n
\choose 2k_n} \prod^n_{j=1} \int {f_j}^{2k_j} \\
&= \sum_{k_1+\ldots+k_n=k} {2k \choose 2k_1} \ldots {2k_n \choose 2k_n}
\prod^n_{j=1} {\|f_j\|_{2k_j}}^{2k_j}.
\end{align*}
Here, equality $(*)$ holds by independence of $f_j$'s and $(**)$ holds by
symmetry of $f_j$'s, and we use a convention that ${\|f_j\|_0}^0 = 1$.
\end{pf}
Notice that when $p=4$, Lemma~\ref{sum} implies:

\begin{cor} \lb{p4}
The uncomplemented subspace $Y_4$ of $L_4$ built in \cite[Corollary after
Proposition~5]{R70} is isometric to a certain complemented subspace of
$L_4$ spanned by $3$-valued independent symmetric random variables.
\end{cor}

\begin{pf}
Let $g_1, \ g_2, \ \ldots$ be a sequence of independent symmetric $\{+1, \
-1, \ 0\}$-valued random variables defined by Rosenthal in \cite[Corollary
after Proposition~5]{R70}, that is, $\int |g_{2n-2}|=1$ and $\int
|g_{2n-3}|=(n \log^2n)^{-1}$ for all $n \geq 2$.  Let
$f_{n-1}=g_{2n-3}+n^{-\frac 12}g_{2n-2}$ and $Y_4= \span \{f_{n-1}\} \subset
L_4$.

Then, by Lemma~\ref{sum},
\begin{align*}
{\|f_{n-1}\|_2}^2 &= {\|g_{2n-3}\|_2}^2 + \|n^{-\frac 12} g_{2n-2}\|_2^2 =
\frac 1{n \log^2n} + \frac 1n\\
\intertext{and}
{\|f_{n-1}\|_4}^4 &= {\|g_{2n-3}\|_4}^4 + 2{\|g_{2n-3}\|_2}^2 \
{\|n^{-\frac 12}g_{2n-2}\|_2}^2 + {\|n^{-\frac 12}g_{2n-2}\|_4}^4 =\\
&= \frac 1{n \log^2n}+2 \cdot \frac 1{n \log^2n} \cdot \frac 1n + \frac
1{n^2}\\
\intertext{Define}
a_n &= \bigg( \frac {n+2+ \log^2n}{n(1+ \log^2n)} \bigg)^{\frac 12}\\
\nu_n &= \frac{1+2 \log^2n+ \log^4n}{n \log^2n+2 \log^2n+ \log^4n} \qquad
\text{ for } n \geq 2.
\end{align*}

Let $h_n$ be $\{1, \ 0, \ -1\}$-valued symmetric independent random
variables with $\int |h_{n-1}|=\nu_n$.

Then
\begin{align*}
{\|a_nh_{n-1}\|_2}^2 &= {a_n}^2 \ \nu_n = {\|f_{n-1}\|_2}^2\\
{\|a_n h_{n-1}\|_4}^4 &= {a_n}^4 \ \nu_n = {\|f_{n-1}\|_4}^4.
\end{align*}

Thus, by Lemma~\ref{sum}, $Z_4= \span \{a_n h_{n-1}\} \subset L_4$ is
isometric with $Y_4 = \span \{f_{n-1}\} \subset L_4$.

By \cite[Theorem~4]{R70}, $Z_4$ is complemented in $L_4$, and by
\cite[Corollary after Proposition~5]{R70}, $Y_4$ is not complemented in
$L_4$.
\end{pf}

When $p$ is an even integer larger than 4, our construction becomes more
complicated, but it is also modeled on \cite[Corollary after
Proposition~5]{R70}.

We will construct subspaces $U_p, \ V_p \subset L_p$ as spans of certain
sequences of independent symmetric random variables.  $V_p$ will not be
spanned by $3$-valued random variables as in the case of $p=4$, but it will
be complemented in $L_p$ by a slight modification of the argument in
\cite[Theorem~4]{R70}.  $U_p$ will be uncomplemented and very similar to
the space $Y_p$ from \cite[Corollary after Proposition~5]{R70}.

We will show that $U_p$ and $V_p$ are isometric by showing that the
$2m$-norms of generators of $U_p$ and $V_p$ are the same for all $m=1,
\ \ldots, \ k$, and by applying Lemma~\ref{sum}.

Consider independent random variables $\{g_{j, i}\}^{k=1, \ \infty}_{i=1, \
j=1}$ which are symmetric $\{+1, \ -1, \ 0\}$-valued, and such that $\int
|g_{j,  i}|$ does not depend on $j$ when $1 \leq i \leq k$, say 
$\int |g_{j,i}|= \mu_i$ for $i=1, \ \ldots, \ k$ and denote 
$\int |g_{j,  k+1}|=\nu_j$.


Define
$h_j = \sum\limits_{i=1}^k g_{j,i}$.


\begin{lem}\lb{normh}

For $m, 1 \leq m \leq k, j\in\bbN, $  $ \Vert h_j\Vert^{2m}_{2m} =
H_m(\mu_1,\dotsc,\mu_k)$,
where $H_m$ is defined by 

$$H_m (\mu_1,\dotsc,\mu_k) = \sum^m_{\alpha = 1} \ \
\sum\Sb
S\subset\{1,\dotsc,k\}  \\
{\card}S = \alpha \endSb
\left[C_{m,\al} \left(\prod\limits_{i\in
s}\mu_i\right)\right],$$ 
where 
$C_{m,\alpha} = \sum\limits_{m_1+\ldots+m_\alpha=m-\alpha}
{2m \choose 2(m_1+1)}\ \
{2(m-(m_1+1)) \choose 2(m_2 + 1)}
 \cdots
{2(m_\alpha+1) \choose 2(m_\alpha +1)}$
is a constant independent of $\mu_1,\dotsc ,\mu_k$.
\end{lem}

\begin{pf}
By Lemma~\ref{sum}, $\Vert h_j\Vert_{2m}$ depends only on
$\left(\Vert g_{ji}\Vert^{2l}_{2l}\right)^{\ k\ \ \ \ m}_{i=1, \ l=1}$.
Moreover by \eqref{sum1} we have, since $\Vert g_{ji}\Vert^{2l}_{2l} =
\mu_i$ whenever $l > 0$,


\begin{align*}
\Vert h_j \Vert^{2m}_{2m}
&= \sum_{m_1+\ldots+m_k=m}
{2m\choose 2m_1}
{2(m-m_1) \choose 2m_2}
\cdots
{2m_k \choose 2m_k}
\prod^k_{i=1} \ \Vert g_{ji}\Vert^{2m_i}_{2m_i}\\
&= \sum_{m_1+\ldots+m_k=m}
{2m\choose 2m_1}
\cdots
{2m_k \choose 2m_k}
\prod^k_{i=1} \mu_i^{{\sgn}(m_i)}\\
&= \sum^m_{\alpha=1}
\ \sum\Sb S \subset \{1,\dotsc,k\}\\
          {\card}S = \alpha\endSb
\left(\prod_{i\in S} \mu_i\right)
\cdot
\left[\sum_{m_1+\ldots+m_\alpha=m-\alpha}
{2m \choose 2(m_1+1)}
{2(m-(m_1+1))\choose 2(m_2+1)}
\ldots
{2m_\alpha \choose2m_\alpha}\right]\\
&= \sum^m_{\alpha=1}
\ \sum\Sb S \subset \{1,\dotsc,k\}\\
          {\card}S=\alpha\endSb
          C_{m,\alpha} \left(\prod_{i\in S}\mu_i\right)\\
&= H_m (\mu_1,\dotsc,\mu_k)\ .
\end{align*}
\end{pf}

Next, define
$f_j = h_j + jg_{j,{k+1}}$

\begin{lem}\label{normf}
For $m, 1 \leq m \leq k\ , j\in \bbN $ we have 
$\Vert f_j\Vert_{2m}^{2m} = 
F_m^{(j)}(\mu_1,\dotsc,
\mu_k, \nu_j)$,
where $F_m^{(j)}$ is defined by
$$
F_m^{(j)}(\mu_1,\dotsc, \mu_k,\nu_j) = H_m(\mu_1,\dotsc,\mu_k) +
\nu_j\sum^m_{l=1}
{2m\choose 2l}j^{2l} H_{m-l}(\mu_1,\dotsc,\mu_k)
$$
\end{lem}
\begin{pf}  We again use \eqref{sum1} to get:
\begin{align*}
\Vert f_j \Vert^{2m}_{2m} &= \Vert h_j + jg_{j,{k+1}}\Vert^{2m}_{2m}\\
&= \sum_{m_1+m_2=2m}
{2m\choose 2m_1}
{2m_2\choose 2m_2}
\Vert h_j \Vert^{2m_1}_{2m_1}
\Vert jg_{j,{k+1}}\Vert^{2m_2}_{2m_2}\\
&= \Vert h_j\Vert^{2m}_{2m} + \sum^m_{l=1}
{2m \choose 2(m-l)}
{2l \choose 2l}
\Vert h_j \Vert^{2(m-l)}_{2(m-l)}\ \cdot j^{2l} \Vert
g_{j,{k+1}}\Vert^{2l}_{2l}\\
&= \Vert h_j\Vert ^{2m}_{2m} + \nu_j \sum^m_{l=1}
{2m\choose 2l}j^{2l}\Vert h_j \Vert^{2(m-l)}_{2(m-l)}\\
&= F_m^{(j)}(\mu_1,\dotsc,\mu_k, \nu_j)\ \ .
\end{align*}
\end{pf}

Our next goal is to show that it is possible to choose
$\mu_1,\dotsc,\mu_k$ and
$\ti\mu_1,\dotsc, \ti\mu_k, \ti\nu_j$
so that
$H_m(\mu_1,\dotsc, \mu_k) = F_m^{(j)}(\ti\mu_1,\dotsc, \ti\mu_k,\ti\nu_j)$
for all
$m = 1,\dotsc,k$, $j\in\bbN$, and, consequently,
so that
$V_p =  {\span}\{h_j\}$
and
$U_p = {\span}\{\ti f_j\}$ are isometric.  

For that, fix $j\in \bbN$ and define
$F^{(j)} : \bbR^{k+1} \lra \bbR^k$ by:
$$
F^{(j)}\left(\mu_1,\dotsc,\mu_k, \nu_j\right) =
\left(F_m^{(j)}(\mu_1,\dotsc,\mu_k,\nu_j)\right)_{m=1}^k\ .
$$

First we will need:

\begin{lem}\lb{rank}   When $\mu_1 > \mu_2 >\ldots>\mu_k> 0$
then for all $j \in \bbN$, $\rank F^{(j)}= k$ at the point
$(\mu_1,\dotsc,\mu_k,0)$.
\end{lem}

\begin{pf}  By Lemma~\ref{normf} for all
$m, 1 \leq m \leq k$ 
\begin{equation}\lb{rank1}
\frac
{\partial F_m^{(j)}}
{\partial\mu_\beta}
\bigg\vert_{(\mu_1,\dotsc,\mu_k,\nu_j)} =
\frac
{\partial H_m}
{\partial {\mu_\beta}}
\bigg\vert_{(\mu_1,\dotsc,\mu_k)} +
\nu_j\sum^m_{l=1}
{2m\choose2l}
j^{2l}\cdot
\frac
{\partial H_{m-l}}
{\partial {\mu_\beta}}
\bigg\vert_{(\mu_1,\dotsc,\mu_k)}\ .
\end{equation}
Thus
$$
\frac
{\partial F_m^{(j)}}
{\partial\mu_\beta}
\bigg\vert_{(\mu_1,\dotsc,\mu_k0)} =
\frac
{\partial H_m}
{\partial {\mu_\beta}}
\bigg\vert_{(\mu_1,\dotsc,\mu_k)}
$$
By Lemma~\ref{normh} we get, that, with a convention
$\prod_{i\in\emp}\mu_i = 1$,

\begin{equation}\lb{partial}
\begin{split}
\frac{\partial H_m}{\partial{\mu_\beta}}\bigg\vert_{(\mu,\dotsc,\mu_k)}
&= \sum^m_{\alpha = 1}\ \sum\Sb S \subset\{1,\dotsc,k\}\\
{\card}S=\alpha\\
\beta\in S \endSb \ \left[C_{m,\alpha}\left(\prod_{i\in S\backslash
\{\beta\}} \mu_i \right)\right]\\
&= \sum^m_{\alpha=1} C_{m,\alpha} \left(\sum\Sb S
\subset\{1,\dotsc,k\}\backslash\{\beta\}\\ {\card}S=\alpha-1\endSb
\left(\prod_{i\in S}\mu_i\right)\right)
\end{split}
\end{equation}

To simplify equation ~\eqref{partial} we will need the following notation:
for any $\alpha,\be : 0 \leq \alpha, \be \leq k$, let

\begin{align*}
P_\alpha &= \sum\Sb S\subset\{1,\dotsc,k\}\\ \text{card}S=\alpha \endSb
\left(\prod_{i\in S}\mu_i\right)\\
P_{\beta,\alpha} &= \sum\Sb S\subset \{1,\dotsc,k\}\backslash \{\beta\}\\
\text{card}S = \alpha \endSb \left(\prod_{i\in S}\mu_i\right)
\end{align*}

Notice that
$P_0 = P_{\beta,0} = 1$, and
\begin{align*}
P_{\beta,\alpha} &= P_\alpha - \mu_\beta \sum\Sb S \subset
\{1,\dotsc,k\}\backslash \{\beta\} \\ \text{card}S = \alpha -1 \endSb
\left(\prod_{i \in S}\mu_i\right) = P_\alpha - \mu_\beta
P_{\beta,\alpha-1}\\
&= P_\alpha - \mu_\beta P_{\alpha-1} + \mu_\beta^2P_{\beta,\alpha -2}\\
&= \sum^\alpha_{t=0} (-1)^t \mu_\beta^t P_{\alpha - t}\ \ .
\end{align*}

Therefore equation~\eqref{partial} becomes:
\begin{align*}
\frac
{\partial H_m}
{\partial{\mu_\beta}}
\bigg\vert_{(\mu_1,\dotsc,\mu_k)}
&= \sum^m_{\alpha = 1} C_{m,\alpha} P_{\beta, \alpha-1}\\
&= \sum^m_{\alpha = 1} C_{m,\alpha} \sum^{\alpha -1}_{t=0} (-1)^t 
\mu_\beta^t P_{\alpha-t-1}\ \ .
\end{align*}

That is,
\begin{align*}
\frac
{\partial H_1}
{\partial{\mu_\beta}}\bigg\vert_{(\mu_1,\dotsc,\mu_k)} \ \ &=\ \ C_{1,1}\\
\frac
{\partial H_2}
{\partial{\mu_\beta}}\bigg\vert_{(\mu_1,\dotsc,\mu_k)}\ \ &=\ \ C_{2,1} +
C_{2,2}\left(P_1 - \mu_\beta\right)\\
\frac
{\partial H_3}
{\partial{\mu_\beta}}\bigg\vert_{(\mu_1,\dotsc,\mu_k)}\ \ &=\ \ C_{3,1} +
C_{3,2}\left(P_2 - \mu_\beta P_1 + \mu_\beta^2\right)\ .
\end{align*}
\noindent
and so on.

Thus the first row of the matrix of the derivative can be reduced to the
row of ones $(1,\dotsc,1)$.  Using the first row we can reduce the 
second row
to the form
$(\mu_1, \mu_2,\dotsc,\mu_k)$.  Similarly, using the first two rows we
reduce the third row to
$(\mu_1^2, \mu_2^2,\dotsc,\mu_k^2)$ and so on.

Thus the matrix of the derivative for $F^{(j)}$ at $(\mu_1,\dotsc,\mu_k)$
reduces
to the Vandermonde matrix
$$
V =
\pmatrix
1     &1     &\cdots     &1\\
\mu_1     &\mu_2     &\cdots    &\mu_k\\
\vdots\\
\mu_1^{k-1}     &\mu_2^{k-1}     &\cdots     &\mu_k^{k-1}
\endpmatrix
$$

Since
$\mu_1 > \mu_2 >\dotsc> \mu_k > 0$, $\det V \not= 0$ (see for example
\cite{HJ}) and therefore Jacobian of $F^{(j)}$ at $(\mu_1,\dotsc,\mu_k)$ is
nonzero, that is $\rank F^{(j)}$ is maximal at the point
$(\mu_1,\dotsc,\mu_k)$.
\end{pf}

Next we want to show that $\rank F^{(j)}$ is maximal in some large enough
open set.

\begin{lem}\lb{nbhd}
There exist $(\overline\mu_1,\dotsc,\overline\mu_k,0) \in \bbR^{k+1}$ and
$\e_0 > 0$ such that for all $j \in \bbN$, $\rank F^{(j)} = k$ at every point
$(\mu_1\dotsc,\mu_k,\nu_j) \in B_\infty
\left((\overline\mu_1,\dotsc,\overline\mu_k, 0), \e_0j^{2-p}\right)$ (here
$B_\infty(x, r)$ denotes the ball with center $x$ and radius $r$ in
$\bbR^{k+1}$ with the $\ell_\infty-\text{norm}$).
\end{lem}

\begin{pf}

Fix some
$(\overline\mu_1,\dotsc, \overline\mu_k)$ with
$\overline\mu_1 > \overline\mu_2 >\dotsc> \overline\mu_k > 0$,
(for example
$\overline\mu_i = \frac ik$ for $i = 1,\dotsc, k)$ and let
$\overline\e > 0$ be such that
$\overline\e < \frac 12 \min \{(\overline\mu_i - \overline\mu_{i+1}), \mu_k
: i = 1,\dotsc,k-1\}$.  Then, by Lemma~\ref{rank},
$$
{\det}\left[
          \left(
                      \frac{\partial F_m^{(j)}}
                           {\partial{\mu_\beta}}
         \bigg\vert_{(\mu_1,\dotsc,\mu_k,0)}
          \right)
^{\ k\ \ \ \ \ k}_{m=1,\beta=1}
          \right] \not= 0
$$
\noindent
for all
$(\mu_1,\dotsc,\mu_k,0) \in B_\infty
\left((\overline\mu_1,\dotsc,\overline\mu_k,0), \overline\e\right)$.

Since the determinant is a continuous function of entries of a matrix there
exists $\e > 0$ such that
$\e \leq \overline\e$ and 
\begin{equation}\lb{epsilon}
\begin{split}
\text{if }
&\max_{1\leq m,\beta\leq k}\bigg\vert a_{m\beta} - \frac{\partial
F_m^{(j)}}{\partial{\mu_\beta}} \bigg\vert_{(\mu_1,\dotsc,\mu_k,0)}\bigg\vert
< \e
\text{ for some }
(\mu_1,\dotsc,\mu_k,0) \in B_\infty
\left((\overline\mu_1,\dotsc,\overline\mu_k,0), \overline\e\right)\ ,\\
&\text{ then }
\det\left((a_{m\beta})^{\ k\ }_{m=1}, ^{\ k\ }_{\beta =1}\right) \not= 0
\end{split}
\end{equation}

Set
$$
M = \sup\left\{
{2m\choose 2l} \bigg\vert \frac{\partial
H_{m-l}}{\partial{\mu_\beta}}\bigg\vert_{(\mu_1,\dotsc,\mu_k)}\ 
:\left(\mu_1,\dotsc, \mu_k,0) \in B_\infty
((\overline\mu_1,\dotsc,\overline\mu_k,0), \overline\e\right), 1 \leq
\beta,m,l \leq k\right\}
$$
and let

\begin{equation} \lb{nbhd2}
\vert \nu_j\vert < \frac \e {M(k-1) j^{2(k-1)}}
\end{equation}

Then, by \eqref{rank1}

$$
\bigg\vert
\frac
{\partial F_m^{(j)}}
{\partial {\mu_\beta}}
\bigg\vert_{(\mu_1,\dotsc,\mu_k,\nu_j)}
\ \ - \ \
\frac
{\partial F_m^{(j)}}
{\partial{\mu_\beta}}
\bigg\vert_{(\mu_1,\dotsc,\mu_k,0)}\bigg\vert
 =
\bigg\vert \nu_j \sum^m_{l=1} {2m\choose 2l} j^{2l} \frac
{\partial H_{m-l}}
{\partial{\mu_\beta}}\bigg\vert_{(\mu_1,\dotsc,\mu_k)}
\bigg\vert  <  \e
$$
when $m = 1,\dotsc, k-1$ and $(\mu_1,\dotsc, \mu_k,0) \in B_\infty
(\overline\mu, \overline\e)$.  Moreover

$$
\bigg\vert\left(
\frac
{\partial F_k^{(j)}}
{\partial{\mu_\beta}}
\bigg\vert_{(\mu_1,\dotsc,\mu_k,\nu_j)}
-
\nu_j
{2k\choose 2k}
j^{2k}
\frac
{\partial H_0}
{\partial{\mu_\beta}}
\bigg\vert_{(\mu_1,\dotsc,\mu_k)}\right)
 -
\frac
{\partial F_k^{(j)}}
{\partial{\mu_\beta}}
\bigg\vert_{(\mu_1,\dotsc,\mu_k,0)}\bigg\vert
$$

$$
=\ \
\bigg\vert \nu_j \sum^{k-1}_{l=1} {2k\choose 2l}j^{2l} \frac
{\partial H_{k-l}}
{\partial{\mu_\beta}}\bigg\vert_{(\mu_1,\dotsc,\mu_k)}
\bigg\vert < \e  .
$$
Thus, if we define for
$\mu = (\mu_1,\dotsc,\mu_k)$, when $m=1,\dotsc,k-1,
\beta = 1,\dotsc,k
$

\begin{equation*}
a_{m\beta}(\mu) = 
\begin{cases} 
{{\disp
\frac{\partial F_m^{(j)}}{\partial{\mu_\beta}}}
\bigg\vert_{(\mu_1,\dotsc,\mu_k,\nu_j)}}\ \
\text{\ \ \ \ \ \ \ \ \ \ \  when}\ \ m =1,\dots,k-1,\  \beta = 1,\dotsc,k\\
\\
{\disp\frac {\partial F_k^{(j)}}{\partial{\mu_\beta}}
\bigg\vert_{(\mu_1,\dotsc,\mu_k,\nu_j)} - \nu_j j^{2k}} \ \ \ \ 
\text{when}\ \ m = k, \ \beta = 1,\dotsc,k
\end{cases}
\end{equation*}
then, by \eqref{epsilon}
$\det((a_{m\beta}(\mu))) \not= 0$.

Further, define
$$
b_{m\beta}(\mu)= \cases a_{m\beta}(\mu), \ \ \text{when}\ \ m =
1,\dotsc, k-1,\  \beta = 1,\dotsc,k\\
\nu_j j^{2k},\ \ \ \  \ \   \text{when}\ \  m=k,\  \beta = 1,\dotsc,k\ .
\endcases
$$
Then, for all
$\beta = 1,\dotsc, k$; we have
\begin{equation*}
\begin{split}
b_{1 \beta}(\mu) = \frac
{\partial F_1^{(j)}}
{\partial{\mu_\beta}}
\bigg\vert_{(\mu_1,\dotsc,\mu_k,\nu_j)}
&=
\frac
{\partial H_1}
{\partial{\mu_\beta}}
\bigg\vert_{(\mu_1,\dotsc,\mu_k)}
+
\nu_j j^2 \frac{\partial H_0}
{\partial{\mu_\beta}}
\bigg\vert_{(\mu_1\dotsc \mu_k)}
\\ &=
C_{1,1} + \nu_j j^2
=
\frac
{C_{1,1} + \nu_j j^2}
{\nu_j j^{2k}} \cdot b_{k\beta}(\mu)\ .
\end{split}
\end{equation*}

Thus $\det(b_{m\beta}(\mu)) = 0$.

Hence we have
\begin{equation}\lb{nbhd1}
\begin{split}
\det
\left(\left(
\frac{\partial F_m^{(j)}}  {\partial{\mu_\beta}}
\bigg\vert_{(\mu_1,\dotsc,\mu_k,\nu_j)}
\right)_{m,\beta}
\right) &=
\det ((a_{m\beta}{(\mu)})) +
\det ((b_{m\beta}(\mu))) \\
&=
\det ((a_{m\beta}(\mu)))
\not= 0\ .
\end{split}
\end{equation}
and \eqref{nbhd1} holds for all
$\mu_1,\dotsc,\mu_k$ with
$(\mu_1,\dotsc, \mu_k,0) \in B_\infty ((\overline
\mu_1,\dotsc,\overline\mu_k,0),\overline\e)$
and all $\nu_j$ satisfying \eqref{nbhd2}, that is \eqref{nbhd1} hold for
all
$(\mu_1,\dotsc,\mu_k,\nu_j) \in B_\infty ((\overline
\mu_1,\dotsc,\overline\mu_k,0), \delta_j)$ where
$\delta_j = \frac {\e}{M(k-1)} j^{2-2k} \DEF {\e_0}{j}^{2-p}$.
\end{pf}

Now we are ready to construct isometric subspaces
$V_p$ and $U_p$ of $L_p$ as indicated before Lemma~\ref{rank}.

\begin{prop}\lb{construction}
There exist
$\overline\mu_1 >\dotsc>\overline\mu_k > 0$ and
$\delta > 0$ such that for all
$j \in \bbN$
there exists $\mu_1^{(j)} >\dotsc> \mu_k^{(j)} > \delta > 0$ and
$\nu_j$ with $\frac {\delta}2 j^{2-p} < \nu_j < \delta j^{2-p}$ so that
$$
F_m^{(j)}\left(\mu_1^{(j)},\dotsc,\mu_k^{(j)}, \nu_j\right) =
H_m\left(\overline\mu_1,\dotsc,\overline\mu_k\right)
$$
for all $m = 1,\dotsc,k$.
\end{prop}

\begin{pf}
Proposition~\ref{construction} is an immediate consequence of
Lemma~\ref{nbhd} and the general theory of multivariable maps whose
derivatives are surjective.  Namely we will use the following:

\medskip

\noindent
{\bf Preimage
Theorem.}{\it (see e.g. \cite[p.21]{DiffTop}).
Let $X,Y$ be manifolds with $\dim X \geq\dim Y$, $f : X \lra Y$ and
$y \in f(X)$ be a value such that the derivative of $f$ at every point 
$x \in
X$ with $f(x) = y$ is surjective.  Then the preimage $f^{-1}(y)$ is a
submanifold of $X$ with $\dim f^{-1}(y) = \dim X - \dim Y$.}

\medskip

Indeed, for each $j \in \bbN$, consider the function
$F^{(j)}\bigg\vert_
{B_\infty
\left(
(\overline\mu_1,\dotsc,\overline\mu_k,0),\e_{0}{j}^{2-p}\right)}$ ,
where $\overline\mu_1,\dotsc,\overline\mu_k,\e_0$ are the numbers
guaranteed by Lemma~\ref{nbhd}.  Then, by Lemma~\ref{nbhd}, the derivative of
$F^{(j)}$ at every point of
$B_\infty\left(\left(\overline\mu_1,\dotsc,\overline\mu_k,0\right),
\e_{0}{j}^{2-p}\right)$
is surjective.

In particular
$y = F^{(j)}(\overline\mu_1,\dotsc,\overline\mu_k,0)$
satisfies the assumptions of the Preimage Theorem, and thus
$f^{-1}(y)$ is a submanifold of
$B_\infty((\overline\mu_1,\dotsc,\overline\mu_k,0),\e_{0}{j}^{2-p})$
with $\dim f^{-1}(y) = (k+1) - k=1 $.  Thus, to satisfy conclusions of
Proposition~\ref{construction}, it is enough to take $\delta = \min(\e_0,
\overline\mu_k - \e_0)$ (recall that for all $m = 1,\dotsc,k$, \ 
$
F_m^{(j)}(\overline\mu_1,\dotsc,\overline\mu_k,0) =
H_m(\overline\mu_1,\dotsc,\overline\mu_k)
$).
\end{pf}

\begin{cor}\lb{def}
Let
$\overline\mu_i, \mu_i^{(j)}, \nu_j,\ \ i = 1,\dotsc,k, j \in \bbN,$
satisfy conditions of Proposition~\ref{construction} and let
$g_{j,i},\ti g_{j,i}$
be such that
$\int\vert g_{j,i}\vert = \mu_i,\ \int\vert\ti g_{j,i}\vert = \mu_i^{(j)},\
\int\vert\ti g_{j,{k+1}}\vert = \nu_j$.
  
Define $V_p = \span\{h_j = \sum^k_{i=1}g_{j,i}\} \subset L_p$
  and
  $U_p = \span \{f_j  = \sum^k_{i=1}\ti g_{j,i} + j\ti g_{j,k+1}\} \subset
  L_p$.  Then $U_p$ and $V_p$ are isometric.
\end{cor}

\begin{pf} Corollary~\ref{def} follows from Lemma~\ref{sum}, 
\ref{normh},
\ref{normf} and Proposition~\ref{construction}.
\end{pf}

Our next goal is to establish that $V_p$ is complemented and $U_p$ is not
complemented in $L_p$

The proof that $V_p$ is complemented in $L_p$ is essentially the same as the
proof that linear $\span$ of independent symmetric 3-valued random variables
is complemented in $L_p$ \cite[Theorem~4]{R70}. The only difference
is that in the case of  3-valued random variables
$\Vert f \Vert_p \cdot \Vert f \Vert_q \leq \Vert f \Vert_2^2$
and in the case of our generating random variables $h_j$ there exists some
larger constant
$C = C(p)$
so that
$\Vert h_j\Vert_p \cdot \Vert h_j \Vert_q
\leq C \cdot \Vert h_j\Vert_2^2$.

For the convenience of the reader we provide the proof (with only one
change from the original).

\begin{lem} \lb{comp} (cf. \cite[Theorem~4]{R70}).
Let $2 < p < \infty$ and let
$h_1,h_2\ldots$ be an infinite sequence of independent symmetric random
variables so that there exists $C$ with
\begin{equation*}
\Vert h_j \Vert_p \cdot \Vert h_j\Vert_q \leq C \Vert h_j \Vert_2^2
\end{equation*}
for all $j$.  Then
$V_p = \span \{h_j\}$
is complemented in $L_p$.
\end{lem}

\begin{pf}  Following Rosenthal, we define $P$ to be the restriction to
$L_p$, of the orthogonal projection from
$L_2$ onto $V_2$, regarding $L_p \subset L_2$.  Thus $P: L_p \lra L_2$
is given by
\begin{equation*}
P(f) = \sum^\infty_{n=1} \left(\int^1_0 f(x) \frac{h_n(x)\Vert
h_n\Vert_p}{\Vert h_n\Vert_2^2}\ dx\right) \frac{h_n}{\Vert h_n\Vert_p}
\end{equation*}
for all $f \in L_p$.  Now fix $f \in L_p$ and put
\begin{equation*}
x_n = \left(\int_0^1f(x) h_n(x) dx\right) \frac{\Vert h_n
\Vert_p}{\Vert h_n\Vert^2_2} \ \ \ \text{for\ all\ } n\ .
\end{equation*}
Then
\begin{equation*}
\left(
\sum \vert x_n\vert^2w_n^2
\right)
^{\frac {1}{2}} =
\left(
\sum \vert x_n\vert^2
\frac{\Vert h_n\Vert_p^2}{\Vert h_n\Vert_2^2}
\right) ^{\frac 1 2} = \Vert Pf\Vert_2 \leq \Vert f \Vert_2 \leq \Vert
f\Vert_p.
\end{equation*}
Now let $n$ be fixed and let
$c_1,\dotsc,c_n$ be such that
$\sum^n_{j=1}\vert c_j\vert^q \leq 1$.  Then
$$\Vert c_j h_j\Vert_q\Vert c_jh_j\Vert_p \leq 
C\vert c_j\vert \Vert h_j\Vert_2^2$$
and \cite[Lemma~2(b) and the remarks following it]{R70} imply that:
\begin{equation*}
\Big\Vert \sum c_j \frac{h_j \Vert h_j\Vert_p}{\Vert h_j\Vert_2^2}
\Big\Vert^q_q
\leq C^q \sum \vert c_j \vert^q \leq C^q.
\end{equation*}
Thus
\begin{equation*}
\big\vert \sum c_j x_j \big\vert =
\bigg\vert \int^1_0 f(x)
\left(\sum c_j \frac{h_j\Vert h_j\Vert_p}{\Vert h_j\Vert_2^2}\right)
dx \bigg\vert
\leq C \Vert f \Vert_p
\end{equation*}
Hence since $n$
and
$c_1,\dotsc,c_n$
were arbitrary satisfying $\sum \vert c_j \vert^q \leq 1$ we have:
\begin{equation*}
 \left(\sum^\infty_{j=1}\vert
x_j\vert^p\right)^{\frac 1p}
\leq C \Vert f \Vert_p\ .
\end{equation*}
Thus
$(x_j)$ belongs to $X_{p,w}$ and
$\Vert (x_j)\Vert \leq C \Vert f \Vert_p$.  Hence by
\cite[Theorem~3 and remarks following it]{R70},
$Pf = \sum x_jf_j$ belongs to
$V_p$
and
$\Vert P_f\Vert_p \leq C\cdot K_p \Vert f \Vert_p$, where $K_p$
is the constant from Rosenthal's inequality \cite[Theorem~3]{R70}.
\end{pf}

To finish the proof of complementability of $V_p$ we need one more:
\begin{lem}\lb{vpl}
When $p$ is an even integer
$(p = 2k)$ then there exists $C = C(p)$ so that
\begin{equation*}
\Vert h_j \Vert_p \cdot \Vert h_j \Vert_q  \leq C \Vert h_j\Vert_2^2
\end{equation*}
for all $j \in \bbN$.
\end{lem}

\begin{pf}
By \cite [Lemma 2 and Remark 1 after it] {R70} we have
\begin{equation*}
\Vert h_j \Vert_q = \Vert \sum^k_{i=1}g_{j,i}\Vert_q \leq \left(\sum^k_{i=1}
\int \vert g_{j,i}\vert^q\right) ^{\frac 1q}=
\left(\sum^k_{i=1}\overline\mu_i\right) ^{\frac 1q}.
\end{equation*} 
By Lemma~\ref{normh}
\begin{equation*}
\Vert h_j\Vert^2_2
= H_1 \left(\overline\mu_1,\dotsc,\overline\mu_k\right)
= C_{1,1} \sum^k_{i=1} \overline\mu_i
= \sum^k_{i=1}\overline\mu_i
\end{equation*}
and
\begin{align*}
\Vert h_j \Vert_p^p
&= H_k (\overline\mu_1,\dotsc,\overline\mu_k)
= \sum^k_{\alpha = 1} \
\sum\Sb S\subset\{1,\dotsc,k\}\\ \card S = \alpha\endSb
\left[C_{k,\alpha} \left(\prod_{i \in S} \overline\mu_i\right)\right] \\
&{\leq{\hspace{-5mm}}\raisebox{-4mm}{\small (*)}} \sum^k_{\alpha = 1}\ \
\sum\Sb S \subset \{1,\dotsc,k\}\\
\card S = \alpha\endSb
C_{k,\alpha} \overline\mu_1 = \overline\mu_1
\cdot \sum ^k_{\alpha = 1} {k \choose \alpha} C_{k,\alpha} \DEF  
\overline\mu_1 \cdot C_k.
\end{align*}
where inequality $(\ast)$ is valid since
$1> \overline\mu_1>\ldots>\overline\mu_k > 0$
and so, for every
$S \subset \{1,\dotsc,k\},$ $\prod_{i\in S} \overline\mu_i \leq
\overline\mu_1$.

Thus we have:
\begin{align*}
\Vert h_j\Vert_q \cdot \Vert h_j\Vert_p
&\leq \left(\sum^k_{i=1} \overline\mu_i\right)^{\frac 1q} \cdot
\overline\mu_1^{\frac 1p} \cdot C_k ^{\frac 1p} \\
&\leq \left(\sum^k_{i=1} \overline\mu_i\right)^{\frac 1q}  \cdot
\left(\sum^k_{i=1} \overline\mu_i\right)^{\frac 1p} \cdot C_k^{\frac 1p}\\
&= \left(\sum \overline\mu_i\right) \cdot C_k^{\frac 1p} = C_k^{\frac 1p}
\cdot \Vert h_j\Vert_2^2\ .
\end{align*}
\end{pf}

Thus from Lemma~\ref{comp} and \ref{vpl} we obtain:
\begin{cor} \lb{vp}

$V_p$ is complemented in $L_p$.
\end{cor}

\begin{rem}
Notice that the constant $C(p)$ in Lemma~\ref {vpl} goes to infinity very
fast when $p$ goes to infinity.  Thus the constant of complementation of
$V_p$ could be very large for large $p$.   We do not know if it is possible
to construct spaces $U_p$ and $V_p$ in such a way that constant of
complementation of $V_p$ is bounded independently of $p$.
\end{rem}

Our last step is to prove that:
\begin{lem} \lb{up}  Assume that $p$ is an even integer and $p\ge 6$.
Then $U_p$ is not complemented in $L_p$.
\end{lem}

\begin{pf}
This is an immediate consequence  of \cite[Theorem~9]{R70} (cf. also 
\cite{E70}).

Indeed we have
\begin{equation*}
\left(\int \vert \ti g_{j,i} \vert\right)^{\frac{p-2}{2p}} = 
\left(\mu_i^{(j)}\right)^{\frac{p-2}{2p}} > \delta^{\frac{p-2}{2p}}
\end{equation*}
 {for} $i = 1,\dotsc,k, \ j
\in \bbN$ and
\begin{equation*}
\sum^\infty_{j=1} \bigg(\big(\int\vert \ti
g_{j,k+1}\vert\big)^{\frac{p-2}{2p}}\bigg)^{\frac{2p}{p-2}} =
\sum^\infty_{j=1} \nu_j < \frac 2 \delta \sum^\infty_{j=1} j^{2-p} < \infty
\end{equation*}
since $p \geq 4$.

Thus, by \cite [Theorem~4 and remarks on page 283]{R70},
$Z_p\DEF \span \{\ti g_{j,i} : i = 1,\dotsc,k+1, j \in \bbN\} \subset L_p$
is isomorphic to $\ell_p \oplus \ell_2$.

Next set $w_j =  1 /({\nu_j^{\frac{1}{p}}j})$.  By
Proposition~\ref{construction}, \ \ 
$\frac\delta2 j^{2-p} < \nu_j < \delta j^{2-p} $.  Thus
\begin{align*}
\left(\frac 1 \delta\right)^{\frac 1p} j^{-\frac 2p}&<w_j<
\left(\frac 2 \delta\right)^{\frac 1p} j^{-\frac 2p}, \\
\lim_{j\to\infty}w_j
&= 0, \\ 
  \sum^\infty_{j=1} w_j^{\frac {2p}{p-2}} &>
\left(\frac1\delta\right)^{\frac 1p} \cdot \sum^\infty_{j=1} 
j^{\frac{-4}{p-2}} \geq \infty\ , 
\end{align*}
 since $p \geq 6$.

Thus, by \cite[Theorem~9]{R70},
$\wt{\wt X}_{p,w} \DEF \span \left\{e_j + w_j 
\sum^k_{i=1}b_{(j-1)k+i}\right\}$
is not a continuous linear image of   $\ell_p \oplus \ell_2 $ (here $(e_j),
(b_j)$ 
denote
unit vector bases in $\ell_p$ and $\ell_2$, resp.)
and we conclude that $U_p$ is
uncomplemented in $Z_p$ and consequently in $L_p$.
\end{pf}

\begin{rem}
In the case when $p = 4$ Lemma~\ref{nbhd} can be extended to say that
rank$F^{(j)} = 2$ at every point
$(\mu_1,\mu_2,\nu_j)$ such that
$\mu_1 > \mu_2 > 0$ and $\nu_j > - \frac{1}{j^2}$.  Therefore also
Proposition~\ref {construction} can be satisfied with $\nu_j$ arbitrarily
large   when $p = 4$; in particular there exists $\nu_j$ with
\begin{equation*}
\frac {\delta}{2} j^{-\frac{3}{2}} < \nu_j < \delta j^{-\frac{3}{2}}
\end{equation*}
satisfying all other conclusions of Proposition~\ref{construction} and
consequently also those of Corollary~\ref{def}.  When $\nu_j$ is of this
order of magnitude it is not difficult to check that $U_4$ is
uncomplemented in $L_4$ (similarly to the proof of Lemma~\ref{up}), however
we felt that it is much simpler to rely on the illustrative
Corollary~\ref{p4} for the case $p = 4$, and we leave the details of the 
above
mentioned computations to the interested reader.
\end{rem}

\bibliographystyle{plain}
\bibliography{tref}

\end{document}